\documentclass[12pt]{amsart}

\usepackage{amsmath}
\usepackage{amsfonts}
\usepackage{amscd}
\usepackage{amssymb}
\usepackage{amsthm}
\usepackage{mathdots}
\usepackage{mathtools}

\usepackage{MnSymbol}
\usepackage[linktocpage=true]{hyperref}
\usepackage{mathrsfs}

\usepackage{tikz,tikz-cd, color}
\usepackage{adjustbox}
\usetikzlibrary{matrix}
\usetikzlibrary{decorations.pathmorphing}

\tikzset{
  symbol/.style={
    draw=none,
    every to/.append style={
      edge node={node [sloped, allow upside down, auto=false]{$#1$}}}
  }
}

\usepackage{stmaryrd}
\usepackage{enumerate}
\usepackage{enumitem}
\usepackage{xcolor}

\makeatletter
\def\@tocline#1#2#3#4#5#6#7{\relax
  \ifnum #1>\c@tocdepth 
  \else
    \par \addpenalty\@secpenalty\addvspace{#2}%
    \begingroup \hyphenpenalty\@M
    \@ifempty{#4}{%
      \@tempdima\csname r@tocindent\number#1\endcsname\relax
    }{%
      \@tempdima#4\relax
    }%
    \parindent\z@ \leftskip#3\relax \advance\leftskip\@tempdima\relax
    \rightskip\@pnumwidth plus4em \parfillskip-\@pnumwidth
    #5\leavevmode\hskip-\@tempdima
      \ifcase #1
       \or\or \hskip 1em \or \hskip 2em \else \hskip 3em \fi%
      #6\nobreak\relax
   \hfill \hbox to\@pnumwidth{\@tocpagenum{#7}}\par
    \nobreak
    \endgroup
  \fi}
\makeatother

\usepackage{cite}

\newtheorem{thm}{Theorem}
\newtheorem*{thm*}{Theorem}
\newtheorem{lem}[thm]{Lemma}
\newtheorem{prop}[thm]{Proposition}

\newtheorem{thm&defn}[thm]{Theorem \& Definition}
\newtheorem{defn}[thm]{Definition}
\newtheorem{exmp}[thm]{Example}

\newtheorem{lemma}[thm]{Lemma}

\theoremstyle{remark}
\newtheorem{rem}[thm]{Remark}
\newtheorem*{rem*}{Remark}

\numberwithin{equation}{section}
\numberwithin{thm}{section}

\newcommand{\bZ}{\mathbb Z}

\newcommand{\bQ}{\mathbb Q}

\newcommand{\bC}{\mathbb C}

\newcommand{\cA}{\mathcal{A}}

\newcommand{\cC}{\mathcal{C}}
\newcommand{\cD}{\mathcal{D}}
\newcommand{\cF}{\mathcal{F}}
\newcommand{\cH}{\mathcal{H}}
\newcommand{\cM}{\mathcal{M}}
\newcommand{\cN}{\mathcal{N}}
\newcommand{\cO}{\mathcal{O}}
\newcommand{\cQ}{\mathcal{Q}}

\newcommand{\cT}{\mathcal{T}}
\newcommand{\cW}{\mathcal{W}}
\newcommand{\cX}{\mathcal{X}}
\newcommand{\cY}{\mathcal{Y}}
\newcommand{\cZ}{\mathcal{Z}}

\newcommand{\db}{{\rm D}^{\rm b}}

\newcommand{\oo}{\mathcal{O}}

\newcommand{\rH}{\mathrm{H}}

\newcommand{\spec}{{\rm Spec\,}}

\newcommand{\rSS}{{\rm SS}}

\newcommand{\coh}{{\rm Coh}}

\newcommand{\Quot}{{\rm Quot}}

\newcommand{\QHusk}{{\rm QHusk}}

\newcommand{\sym}{{\rm Sym}}

\newcommand{\St}{{\rm St}}

\newcommand{\id}{{\rm id}}

\newcommand{\im}{{\rm im\,}}
\newcommand{\coker}{{\rm coker\,}}

\newcommand{\Coh}{{\rm Coh}}
\newcommand{\Hom}{{\rm Hom}}

\newcommand{\Ext}{{\rm Ext}}

\newcommand{\Db}{{\mathrm{D}^{\mathrm{b}}}}

\newcommand{\cohltilt}{\mathrm{Coh}^{\lambda,\#}(X)}

\newcommand{\ignore}[1]{}

\title{Decorated sheaves and morphisms in tilted hearts}
\author{Yinbang Lin}
\address[]{School of Mathematical Sciences, Tongji University, Shanghai, China}
\email{yinbang\textunderscore lin@tongji.edu.cn}
\author{Sz-Sheng Wang}
\address[]{Shing-Tung Yau Center of Southeast University, Southeast University, Nanjing, China}
\email{sswang@seu.edu.cn}
\author{Bingyu Xia}
\address[]{Morningside Center of Mathematics, Chinese Academy of Science, Beijing, China}
\email{xiabingyu@amss.ac.cn}
\date{}


\begin{document}

\begin{abstract}
    We identify limit stable pairs and stable framed sheaves as epimorphisms and monomorphisms, respectively, in tilts of the standard heart, under suitable conditions. We then identify the moduli spaces with the corresponding Quot spaces, obtaining the projectivity of the Quot spaces in these cases. We also prove a formula in a motivic Hall algebra relating the Quot spaces under a tilt.
\end{abstract}
\maketitle
\section{Introduction}

Decorated sheaves are sheaves with extra structure.
The most famous moduli space of decorated sheaves is Grothendieck's Quot scheme. Quite often, the Quot scheme is oversized. Towards some problems, certain variants are more suitable.
Among them, one is the moduli space of stable pairs. It was used by Thaddeus \cite{Tha94} to calculate the Verlinde numbers. It was also used by Pandharipande and Thomas \cite{PanTho09} to study curve counting on Calabi--Yau 3-folds. For a more recent attempt towards strange duality, see \cite{goller2019rankone}.
We will focus on two variants of quotient sheaves: stable pairs and framed sheaves.

Let $X$ be a nonsingular projective variety with a fixed polarization over an algebraically closed field of characteristic 0. Let $E_0$ be a fixed coherent sheaf on $X$.
By an unfortunate abuse of terminology, we call a sheaf equipped with a morphism $\alpha$ $$(E,\alpha\colon E_0\to E)$$
 a {\em pair}.
There is a family of stability conditions (Definition \ref{stable-pair}) defined on pairs.
When the stability polynomial is large, a pair $(E,\alpha)$ is stable if and only if $\alpha$ is generically surjective. In this case, we call the pair {\em limit stable}.
We will identify stable pairs with quotients in the hearts of certain tilts of the standard $t$-structure (Proposition \ref{pair-quotient}).

Limit stable pairs are also known as quotient husks \cite{Kol0805}. There, Koll\'{a}r studied in a relative setting.
Suppose $X\to S$ is a smooth projective morphism.
We define a torsion pair on $\coh(X)$ (Lemma~\ref{relative-torsion-pair}). The tilting $t$-structure, whose heart is denoted as $\coh^\#(X)$, pulls back to the corresponding tilting $t$-structure on each fiber.
A quotient in $\coh^\#(X)$ corresponds to a family of quotient husks, under suitable conditions (\S\ref{moduli-husk}).
The moduli space satisfies the valuative criteria of separatedness and properness. This can be understood as giving a specialization map between moduli spaces.

There is a dual notion: {\em framed sheaves} \cite{HL95}. A framed sheaf is of the form $$(E,\alpha\colon E\to E_0).$$
There is also a family of stability conditions (Definition \ref{framed-sheaf}) on framed sheaves. We will identify stable framed sheaves with monomorphisms in certain tilted hearts (Proposition \ref{framed-sheaf-injection}).

Over an algebraically closed field of characteristic 0, these moduli spaces coincide with the corresponding Quot spaces \cite{BLM+19}, which are generalizations of Quot schemes.
As a consequence, these Quot spaces are projective.
We also point out that torsion pairs where we tilt the standard heart are defined using the corresponding stability polynomials.

Over a Calabu--Yau 3-fold, when $E_0\cong \oo_X$ and $E$ has dimension 1, the equivalence between stable pairs and quotients in a tilted heart has been obtained by Bridgeland \cite{Bridgeland11} and used to derive a wall-crossing formula between the Donaldson--Thomas (DT) invariants and Pandharipande--Thomas (PT) invariants.
As an immediate application of our results, we also derive a formula relating the moduli space of limit stable pairs and the Quot scheme under certain assumptions (Theorem \ref{App_thm}).

Over a nonsingular curve, stable pairs of the form $(E,\alpha\colon \oo_X\to E)$ have been identified as epimorphisms in the corresponding tilted heart in \cite{rota2020quot}. A variant of the result on framed sheaves has also been obtained there.

We organize the paper as follows.
In \S\ref{torsion-pair}, we review some basic notions on $t$-structures.
In \S\ref{moduli-husk}, we compare quotient husks and quotients in a tilted heart.
In \S\ref{pair}, we identify stable pairs and framed sheaves as epimorphisms and monomorphisms in different tilted hearts.
In \S\ref{wall-crossing}, we review the formalism of motivic Hall algebra and derive the formula relating the Quot spaces under a tilt.

\section{$t$-structures and torsion pairs}\label{torsion-pair}
In this section, we review the basic notions of $t$-structures, hearts, torsion pairs, and tilts, and prove an observation important to us.

A \emph{$t$-structure} $\tau$ \cite{BBDG} on a triangulated category $\cD$ is a pair of full subcategories $(\cD^{\leqslant 0}, \cD^{\geqslant 0})$ satisfying the following conditions:
\begin{enumerate}[label=(\roman*)]
    \item $\Hom (F, G) = 0$ for every $F \in \cD^{\leqslant 0}$ and $G \in \cD^{\geqslant 1}$;

    \item every object $E \in \cD$ fits into an exact triangle
    $$
    \tau^{\leqslant 0} E \to E \to \tau^{\geqslant 1} E \to (\tau^{\leqslant 0} E) [1]
    $$
    with $\tau^{\leqslant 0} E \in \cD^{\leqslant 0}$ and $\tau^{\geqslant 1} E \in \cD^{\geqslant 1}$.
\end{enumerate}
Here, we use the notation $\cD^{\leqslant n} \coloneqq \cD^{\leqslant 0} [- n]$ and $\cD^{\geqslant n} \coloneqq \cD^{\geqslant 0} [- n]$, for $n \in \bZ$. The truncation functors $\tau^{\leqslant n}$ and $\tau^{\geqslant n}$ are similarly defined. The full subcategory $\cA \coloneqq \cD^{\leqslant 0} \cap \cD^{\geqslant 0}$ is called the \emph{heart} of the $t$-structure, which is an abelian category. The $t$-structure is \emph{bounded} if $\cD = \cup_{n, m \in \bZ} \cD^{\leqslant n} \cap \cD^{\geqslant m}$. The cohomology objects of an object $E \in \cD$ with respect to the heart $\cA$ are defined by $\rH_{\cA}^{n} (E) \coloneqq (\tau^{\leqslant n} \tau^{\geqslant n} E) [n]$.

Assume that $\cA$ is an abelian category. 
A pair of full subcategories $(\cT, \cF)$ of $\cA$ is called a \emph{torsion pair} if $\Hom (\cT, \cF) = 0$ and every object $E \in \cA$ fits into an exact sequence
\begin{equation}\label{tor_pair}
    0 \to T \to E \to F \to 0
\end{equation}
with $T \in \cT$ and $F \in \cF$. It follows that $\cT$ and $\cF$ are respectively closed under taking quotients and sub-objects. Note that if $(\cT, \cF)$ is a torsion pair, then $\cF = \cT^{\perp}$ is the right orthogonal to $\cT$ in $\cA$.

\begin{exmp}\label{coh<=}
Let $X$ be a noetherian scheme of dimension $n$. For any integer $0 \leqslant d <n$, we define
\[
    \coh_{\leqslant d} (X) \coloneqq \{ E \in \coh (X) \mid \dim E \leqslant d \}.
\]
Here and henceforward, we use $\dim E$ to denote the dimension of the support of the sheaf $E$. When $d = 0$, we simply write $\coh_{0} (X)$. Then $(\cT, \cT^{\perp})$ is a torsion pair in the abelian category $\coh (X)$ (cf.~\cite[Definition 1.1.4]{HL10}).
\end{exmp}

Let $\cA \subset \cD$ be the heart of a $t$-structure and $(\cT, \cF)$ a torsion pair in $\cA$. We can \emph{tilt} the $t$-structure  on $\cD$ to obtain a new one. The resulting $t$-structure has heart $\cA^{\#} = \left<\cF, \cT [-1]\right>$, the extension-closure. Also, the shift $\cA^{\#} [1]$ can be described as
\[
    \cA^{\#} [1] = \left\{ E \in \cD \mid \rH_{\cA}^0 (E) \in \cT,\  \rH_{\cA}^{- 1} (E) \in \cF,\  \rH_{\cA}^{i} (E) = 0 \text{ for } i \neq 0,-1 \right\}.
\]

The following observation is important for us.

\begin{lemma} \label{surjtil}
Let $\cA$ be the heart of a bounded $t-$structure on a triangulated category $\cD$ and $(\cT, \cF)$ a torsion pair of $\cA$. Let $\cA^{\#}$ be the tilt in $\cA$ at this pair.
\begin{enumerate}[label=(\roman*)]
    \item\label{surjtil-r} If $E_0 \in \cF$, then a morphism $\alpha \colon E_0 \to E$ in $\cA^{\#}$ is an epimorphism if and only if $\alpha$ is a morphism in $\cA$ with $E \in \cF$ and the cokernel $\coker (\alpha)$ taken in $\cA$ lies in $\cT$.

    \item\label{surjtil-l} If $E_0 \in \cT$, then a morphism $\beta \colon E \to E_0$ in $\cA^{\#}[1]$ is a monomorphism if and only if $\beta$ is a morphism in $\cA$ with $E \in \cT$ and the kernel $\ker(\beta)$ taken in $\cA$ lies in $\cF$.
\end{enumerate}
\end{lemma}

\begin{proof}
We only show \ref{surjtil-l}, as the proof of \ref{surjtil-r} is similar (cf.~\cite[Lemma 2.3]{Bridgeland11}). Given $E_0 \in \cT$ and a short exact sequence in $\cA^{\#}[1]$
\begin{equation*}
    0 \to E \xrightarrow{\beta} E_0 \to Q \to 0,
\end{equation*}
we have
\begin{equation} \label{longexc-ltil}
    0 \to \rH^{-1}_{\cA} (E) \to 0 \to \rH^{-1}_{\cA} (Q) \to \rH^0_{\cA} (E) \xrightarrow{\beta} E_0 \to \rH^0_{\cA} (Q) \to 0
\end{equation}
by taking cohomology with respect to the original $t$-structure. Note that an object $F \in \cD$ lies in $\cA^{\#}[1] \subset \cD$ precisely when $\rH^{-1}_{\cA} (F) \in \cF$, $\rH^{0}_{\cA} (F) \in \cT$ and $\rH^{i}_{\cA} (F) = 0$ for $i \neq 0, -1$. It follows that $\rH^i_{\cA} (E) = 0$ for all $i \neq 0$, and thus $E \in \cA \cap \cA^{\#}[1] = \cT$.
Moreover, $\ker (\beta) = \rH^{-1}_{\cA} (Q) \in \cF$.

For the converse, take a morphism $\beta \colon E \to E_0$ in $\cA$ with $E \in \cT$, $\ker (\beta) \in \cF$ and embed it in a distinguished triangle
\begin{equation}\label{tri-ltil}
    E \xrightarrow{\beta} E_0 \to Q \to E [1].
\end{equation}
Since $\cT$ is closed under taking quotients and $E_0 \in \cT$, the same long exact sequence \eqref{longexc-ltil} then shows that $\rH^0_{\cA} (Q) \in \cT$, $\rH^{- 1}_{\cA} (Q) = \ker (\beta) \in \cF$ and $\rH^{i}_{\cA} (Q) = 0$ for $i \neq 0,-1$. Thus $Q \in \cA^{\#}[1]$. It follows that \eqref{tri-ltil} defines a short exact sequence in $\cA^{\#}[1]$, and hence $\beta$ is a monomorphism in $\cA^{\#}[1]$.
\end{proof}

\section{Quot spaces and moduli of quotient husks}\label{moduli-husk}
We construct a family of $t$-structures for a smooth projective morphism, and compare the Quot space in the sense of \cite[\S 11]{BLM+19} with the moduli space of quotient husks constructed by Koll\'{a}r \cite{Kol0805}.

Over a field $k$, let $S$ be a regular noetherian scheme of finite type and \[f\colon X\rightarrow S\] be a smooth projective morphism, with a fixed relatively ample line bundle.
On each fiber $X_{s}$ for $s\in S$, we can define the following tilting $t$-structure for some fixed non-negative integer $m$:
\begin{align*}
    \mathcal{T}_{s} &= \{E\in \coh(X_{s}) \mid    \dim  E\leqslant m\},\\
    \mathcal{F}_{s} &= \{F\in \coh(X_{s}) \mid \Hom(E,F)=0,\forall E\in \mathcal{T}_{s}\}.
\end{align*}
Then $(\mathcal{T}_{s},\mathcal{F}_{s})$ defines a torsion pair of $\coh(X_{s})$. If we denote the standard $t$-structure on $\db(X_{s})$ as $(D^{\leqslant0}_{s},D^{\geqslant0}_{s})$, the tilting $t$-structure is defined by
\begin{align}\label{fiber-t-str}
D^{\#, \leqslant0}_{s} = \{E\in D^{\leqslant1}_{s} \mid \rH^{1}(E)\in \mathcal{T}_{s}\}\mbox{ and }
D^{\#, \geqslant0}_{s} = \{E\in D^{\geqslant0}_{s} \mid \rH^{0}(E)\in \mathcal{F}_{s}\}.
\end{align}
The cohomologies are taken with respect to the standard $t$-structure. The heart of the $t$-structure is $\mathcal{A}_{s}^{\#}=D^{\leqslant0}_{s}\cap D^{\geqslant0}_{s}$. It is an abelian category, so we can talk about epimorphisms in it.

We will show that the family of $t$-structures defined above integrates over $S$, namely, they are pullbacks of a $t$-structure on the total space $X$.

We have two full subcategories of $\Db(X)$:
\begin{align*}
\mathcal{T} &= \{E\in \coh(X) \mid \forall s\in S, \dim E_s\leqslant m\},\\
\mathcal{F} &= \{F\in \coh(X) \mid \Hom(E,F)=0, \forall E\in \mathcal{T}\}.
\end{align*}
For brevity, we say objects in $\mathcal{T}$ have {\em relative dimension $\leqslant m$ over }$S$.
\begin{lem}\label{relative-torsion-pair}
$(\mathcal{T},\mathcal{F})$ is a torsion pair of $\coh(X)$.
\end{lem}
\begin{proof}
This is similar to the absolute case. For any coherent sheaf $E$, We can take its maximal subsheaf $E'$ whose support has a relative dimension $\leqslant m$ over $S$. Then the quotient $E/E'$ belongs to $\cF$.
\end{proof}
The torsion pair $(\mathcal{T}, \mathcal{F})$ induces the following $t$-structure on $\db(X)$:
\begin{align}\label{nonconst-t-str}
D^{\#, \leqslant0} = \{E\in D^{\leqslant1} \mid \rH^{1}(E)\in \mathcal{T}\}\mbox{ and }
D^{\#, \geqslant0} = \{E\in D^{\geqslant0} \mid \rH^{0}(E)\in \mathcal{F}\}
\end{align}
where $(D^{\leqslant0},D^{\geqslant0})$ is the standard $t$-structure on $\db(X)$.
By \cite[Theorem 5.7]{BLM+19}, there is a canonical way to pullback this $t$-structure to $\db(X_{T})$ for an arbitrary morphism $T\rightarrow S$, which is {\em faithful} \cite{Kuz11} with respect to $f$ since $f$ is flat.
Here, \[X_T=T\times_S X.\]Suppose $s$ is a point of $S$ (possibly nonclosed) and we denote
\begin{equation*}
((D^{\#})_{s}^{ \leqslant0},(D^{\#})_{s}^{ \geqslant0})\end{equation*}
to be the pullback $t$-structure on the fiber $X_{s}$. We have the following:
\begin{lem}Pullbacks of the $t$-structure agree with the fiberwise defined ones, namely,
$(D^{\#})_{s}^{ \leqslant0}=D^{\#, \leqslant0}_{s}$ and $(D^{\#})_{s}^{ \geqslant0}=D_{s}^{\#, \geqslant0}$.
\end{lem}

\begin{proof}
If $s$ is a closed point, we can use \cite[Theorem 5.7 (3)]{BLM+19}. Let $T=\{s\}$ and $\phi$ be the closed embedding of $\{s\}$ into $S$. Take $a=-\infty$, $b=0$, then 
$D^{(-\infty,0]}_{s}=(D^{\#})_{s}^{ \leqslant0}$ and $D^{(-\infty,0]}=D^{\#,\leqslant0}$. Let $i_{s} \colon X_{s}\rightarrow X$ be the inclusion of the fiber $X_{s}$ in $X$, we have
\begin{align*}
       (D^{\#})_{s}^{ \leqslant0}&=\{F\in\Db(X_{s})\mid i_{s*}F\in D^{\#, \leqslant0}\}\\
       &=\{F\in\Db(X_{s})\mid F\in D^{\leqslant1}, \rH^{1}(F)\in \mathcal{T}\}=D^{\#, \leqslant0}_{s}.
\end{align*}
The other equality can be shown similarly.

For a non-closed point $s$, we can still apply \cite[Theorem 5.3 \& 5.7]{BLM+19} and obtain the result.
\end{proof}

\begin{rem}
Suppose $X=S\times Y\to S$ is a trivial family with $Y$ over an algebraic closed field $k$, the $t$-structure (\ref{nonconst-t-str}) may {\it not} be the pullback of the one on $\db(Y)$, which is similar to (\ref{fiber-t-str}). Namely, it is not constant in the sense of \cite{AP06,Pol07}, as illustrated in the following example.
\end{rem}
\begin{exmp}
Let $Y=\mathbb{P}^1$ and $S=\mathbb{P}^1$. Let $\Delta\cong \mathbb{P}^1\subset X=S\times Y$ be the diagonal. Then $\oo_\Delta[-1]$ lies in the the heart $D^{\#,\leqslant 0}\cap D^{\#,\geqslant 0}$. However, it does not lie in the heart of the pullback of the corresponding $t$-structure, according to \cite[Lemma 3.3.2]{Pol07}.
\end{exmp}

Let $F_{0}\in\mathcal{F}\subseteq\coh(X)$ be pure and flat over $S$.\footnote{The Quot functor can be defined in a much broader context, but let us restrict to this situation, which is enough for our purpose.}
Let $T$ be an $S$-scheme and $F_{0T}$ denote the pullback of $F_0$ to $X_T$.
A family of epimorphisms in the tilted hearts of the family of $t$-structures (\ref{fiber-t-str}), parametrized by $T$, is a morphism $F_{0T}\rightarrow F$ in $\db(X_{T})$ satisfying the following conditions.
\begin{enumerate}[label=(\alph*)]
    \item $F$ is flat over $T$, i.e. the derived pullback $F_{t}$ of $F$ to a fiber $X_{t}$ lies in the tilted heart $\coh^{\#}(X_{t})$ of the induced $t$-structure on $\Db(X_{t})$.
    \item The morphism $F_{0t}\rightarrow F_{t}$ in $\coh^{\#}(X_{t})$ is an epimorphism for all $t\in T$.
\end{enumerate}
Here,
\[X_t=\spec \kappa(t)\times_S X.\]
\begin{rem}
In this setting, $F$ is automatically a sheaf because its derived pullback $F_{t}$ to each fiber is a quotient of  $F_{0t}$ in $\coh^{\#}(X_{t})$ and by Lemma \ref{surjtil}, it must be a sheaf in $\mathcal{F}_{t}$ .
\end{rem}

Let
\begin{equation*}
    \mathrm{Quot}_{f}^\#(F_{0},P)\colon (S\mbox{-Schemes})^{\mathrm{op}}\rightarrow\mathrm{Sets}
\end{equation*}
be the functor that sends an $S$-scheme $T$ to the set of families of epimorphisms in the tilted hearts where the quotient has Hilbert polynomial $P$. If the family of $t$-structures universally satisfies openness of flatness, then the functor is represented by an algebraic space locally of finite presentation over $S$ \cite[Proposition 11.6]{BLM+19}.

We next recall the notion of a {\em husk}. A husk of a coherent sheaf $F$ on a normal scheme $Y$ is a homomorphism $q\colon F \rightarrow E$ such that
\begin{enumerate}[label=(\alph*)]
    \item $q$ is an isomorphism on all $n$-dimensional points, where $n=\dim F$;
    \item $E$ is pure of dimension $n$.
\end{enumerate}
A {\em quotient husk} of a fixed coherent sheaf $E_0$ is a homomorphism $q\colon E_0\rightarrow E$ such that it factors as $E_0\rightarrow F\rightarrow E$ where the first arrow is an epimorphism and the second arrow is a husk.
For a morphism $f\colon X\rightarrow S$, one can similarly define the notion of a family of quotient husks to be a morphism $q\colon E_0\rightarrow E$ such that
\begin{enumerate}[label=(\alph*)]
    \item $E$ is pure and flat over $S$;

    \item On each fiber $X_{s}$, the homomorphism $q_{s}\colon E_{0s}\rightarrow E_{s}$ is a quotient husk.
\end{enumerate}
We say a sheaf $E$ on $X$ is {\em pure over} $S$ if for every $s\in S$, the restriction $E_s$ is pure of the same dimension. It is immediate to check that this definition of a family of quotient husks has base-change property, therefore we can define a moduli functor of families of quotient husks
\[\mathrm{QHusk}_{f}(E_0, P)\colon (S\mbox{-Schemes})^{\rm op}\to \mbox{Sets}\]
where $P$ is a fixed polynomial with rational coefficients and $E_0\in\coh(X)$.
It sends an $S$-scheme $T\rightarrow S$ to the set of families of quotient husks $f_{T}^{*}E_0\rightarrow E_T$ 
on $X_{T}$ such that when restricted to each fiber $X_t$ for $t\in T$, the
Hilbert polynomial of $E_t$ is $P$.
Then we have the following existence theorem on the moduli space of quotient husks by Koll\'{a}r:
\begin{thm}
$\mathrm{QHusk}_{f}(E_0, P)$ is represented by a proper 
algebraic space.
\end{thm}

The following proposition shows that the notion of a family of quotient husks is equivalent to a family of epimorphisms with respect to the family of $t$-structures (\ref{fiber-t-str}).
\begin{prop} \label{husk-quotient}
Assume $E_{0}, E$ are two sheaves flat over $S$ and the restriction $E_s$ has $(m+1)$-dimensional support, then a homomorphism of sheaves $q\colon E_0\rightarrow E$ is a family of quotient husks if and only if it is a family of quotients with respect to the family of $t$-structures (\ref{fiber-t-str}).
\end{prop}

\begin{proof}
Suppose $q\colon E_{0}\rightarrow E$ is a family of quotient husk, then for every point $s\in S$, the restriction $q_{s}\colon E_{0s}\rightarrow E_{s}$ is a quotient husk on $X_{s}$. In particular, $E_{s}$ is a pure sheaf with $(m+1)$-dimensional support. This implies $E_{s}$ is contained in $\mathcal{F}_{s}$. The fact that $q_{s}$ is surjective at all $(m+1)$-dimensional points implies its cokernel is supported in a locus of dimension $m$ or less. Together with Lemma \ref{surjtil}, this shows that $q_{s}$ is an epimorphism in $\coh^{\#}(X_{s})$.

Conversely, if $q$ is a family of quotients, then on each fiber, the homomorphism $q_{s}$ factors as $E_{0s}\rightarrow\mathrm{Im}(q_{s})\rightarrow E_{s}$ where the first arrow is an epimorphism. By Lemma \ref{surjtil}, the cokernel of $q_{s}$ is in $\mathcal{T}_{s}$ so it is supported in a locus of dimension $\leqslant m$. Therefore the second arrow above has to be an isomorphism at all $(m+1)$-dimensional points. Also by Lemma \ref{surjtil}, $E_{s}$ is contained in $\mathcal{F}_{s}$ and therefore does not contain any subsheaf supported on a locus of dimension $m$ or less. This implies $E_{s}$ is pure and $E$ is pure over $S$. This finishes the proof that $q_{s}$ is a quotient husk.
\end{proof}

Let $E_0=F_0$ and $P$ be a polynomial of degree $m+1$, then the functor $\Quot_f^\#(E_0, P)$ is the same as the functor $\mathrm{QHusk}_f(E_0, P)$, that is, we have the following proposition.
\begin{prop}\label{husk-quotientCor}
Given an $S$-scheme $T$, $\Quot_f^\#(E_0, P)(T)=\mathrm{QHusk}_f(E_0, P)(T)$.
\end{prop}

\begin{proof}
Suppose we are given an element $\alpha_T\colon E_{0T}\to E$ in $\Quot_f^\#(E_0,P)(T)$. For $t\in T$, the restriction $\alpha_t\colon E_{0t}\to E_t$ lies in the heart $\coh^{\#}(X_t)$ of the pullback $t$-structure. 
By Lemma \ref{surjtil}, $E_t\in \coh(X_t)$. Thus, $E$ is a coherent sheaf flat over $T$. On the other hand, $\alpha_t$ has cokernel (taken in $\coh( X_t)$) in $\mathcal{T}_{t}$. Therefore, $\alpha_T$ is a family of quotient husks.

Given a family of quotient husks $\alpha_T\colon E_{0T}\to E$, the restriction $\alpha_t$ is a quotient husk. By Proposition~\ref{husk-quotient} and Lemma~\ref{surjtil}, $E$ is flat over $T$ and $\alpha_t$ is a quotient in $\coh^\#( X_t)$.
\end{proof}

By identifying the moduli functors, we can easily obtain the projectivity of the Quot space in some cases.

\begin{thm}
Let $k$ be an algebraic closed field of characteristic 0 and $f\colon X\to S=\spec k$ be a nonsingular projective scheme. Let $P$ be a polynomial of degree $m+1$. Then the Quot space $\Quot_f^\#(E_0,P)$ is projective.
\end{thm}
\begin{proof}
This follows from the projectivity of $\mathrm{QHusk}_f(E_0, P)$ in this set-up, which is obtained via a geometric invariant theoretic construction \cite{Lin18}.
\end{proof}

\begin{rem}
Quotient husks are also known as limit stable pairs \cite{Lin18}. Assuming the universal openness of flatness for the family of $t$-structures, one would be able to obtain the projectivity of the Quot space over a general base $S$ in characteristic 0, by carrying out a GIT construction \cite[Remark 4.6]{Lin18}.
\end{rem}


For stable pairs with respect to a smaller stability condition, we can also identify them as quotients in the heart of a tilting $t$-structure. We will study this next.

\section{Stable pairs and framed sheaves}\label{pair}

Over an algebraically closed field of characteristic 0, let $X$ be a nonsingular projective variety with a fixed polarization $\mathcal{O}_X(1)$. Let $E_0\in \Coh(X)$ be fixed. Let $P\in \mathbb{Q}[m]$ be a fixed polynomial of degree $d$, which will be used as a Hilbert polynomial.
Let $\delta,\tau\in \mathbb{Q}[m]$ be polynomials with positive leading coefficients.

Given a coherent sheaf $E$ on $(X,\mathcal{O}_X(1))$, we denote its Hilbert polynomial by $P_E$, its multiplicity by $r_E$ and its reduced Hilbert polynomial by $p_E=P_E/r_E$.

\subsection{Stable pairs}
We consider homomorphisms of the form $E_0\to E$.
\begin{defn}\label{stable-pair}
A {\em pair}
$$(E, \alpha\colon E_0\to E)$$ with $\alpha\neq 0$
is $\delta$-{\em stable} if $E$ is pure and for every subsheaf $F\subset E$,
\begin{enumerate}[label=(\roman*)]
    \item $p_F+\delta/r_F< p_E+\delta/r_E$ if $\im \alpha\subset F$,
    \item $p_F< p_E+\delta/r_E$ otherwise.
\end{enumerate}
\end{defn}
We can replace the strong inequalities by weak inequalities to define $\delta$-semistability.
Stability can be equivalently defined in terms of quotients.
When $\deg \delta\geqslant \dim E$, a $\delta$-stable pair is called a {\em limit stable pair}, which is the same as a quotient husk.

\begin{thm}[Y. Lin
]
There is a projective coarse moduli space $S_{E_0}(P,\delta)$ of S-equivalence classes of $\delta$-semistable pairs with Hilbert polynomial $P$. It contains an open subscheme $S_{E_0}^{\rm s}(P,\delta)$ as the fine moduli space of $\delta$-stable pairs.
\end{thm}

We consider $\delta$-stable pairs $(E,\alpha)$ with $P_E=P$ for a small $\delta$:
$$\deg \delta<\deg P=d.$$
Let $r=r_E$ and
\begin{equation*}
    \lambda=\frac{P+\delta}{r}.
\end{equation*}
We define a torsion pair $(\mathcal{T}^\lambda, \mathcal{F}^\lambda)$ on $\Coh(X)$:
\begin{eqnarray*}
\mathcal{T}^\lambda &=&\{E\in \Coh(X) \mid \dim E\leqslant d\mbox{ and } \forall \mbox{ quotient sheaf }E\twoheadrightarrow G,\ \dim G<d\mbox{ or }\ p_G>\lambda\} \quad \mbox{and}\\
\mathcal{F}^\lambda &=&\{E\in \Coh(X) \mid \forall \mbox{ subsheaf }F\subset E\mbox{ with }\dim F\leqslant d, \ F\mbox{ is pure of dim. }d\mbox{ and }\ p_F\leqslant\lambda\}.
\end{eqnarray*}
Then, we denote the heart of the tilting $t$-structure as $\Coh^{\lambda,\# }(X)$.

\begin{prop}\label{pair-quotient}
 Suppose $E$ has Hilbert polynomial $P$ and $E_0\in \mathcal{F}^\lambda$.
 We also assume that $\delta$ is not a critical value, namely, there are no strictly semistable pairs with respect to $\delta$. Then, $\alpha\colon E_0\to E$ is an epimorphism in $\cohltilt$ if and only if $(E,\alpha\colon E_0\to E )$ is a $\delta$-stable pair.
\end{prop}

\begin{proof}
Suppose that $\alpha$ is an epimorphism in $\cohltilt$. By the assumption $E_0 \in \cF^{\lambda}$ and Lemma \ref{surjtil}, we know that $\alpha$ is a morphism in $\coh (X)$ with $E \in \cF^{\lambda}$ and $\coker (\alpha) \in \cT^{\lambda}$.  Given a quotient $q\colon E\twoheadrightarrow G$ in $\Coh(X)$, if $q\circ \alpha=0$, then $G$ is a quotient of $\coker (\alpha)$. Thus, $p_G>\lambda$. If $q\circ \alpha\neq 0$, let $F=\ker q$. Then $p_F\leqslant \lambda$. Since $\delta$ is not critical, this is a strict inequality. Therefore, $\lambda<(P_G+\delta)/r_G$. We have shown $(E,\alpha)$ is a $\delta$-stable pair.

Conversely, suppose $(E,\alpha\colon E_0\to E)$ is a $\delta$-stable pair. Then for every subsheaf $F\subset E$, $p_F<\lambda$. Thus, $E\in \cF^{\lambda}$. On the other hand, $\coker (\alpha)$ has dimension $\leqslant d$. Given a quotient $G$ of $\coker (\alpha)$, it is also a quotient of $E$ and the composition $E_0\to E\to G$ is zero. Thus, $p_G>\lambda$. Hence, $\coker (\alpha)\in \mathcal{T}^\lambda$. Again by Lemma \ref{surjtil}, $\alpha$ is an epimorphism in $\cohltilt$.
\end{proof}

Similar to the quotient husks case, by identifying two moduli functors, we can identify the moduli space with the corresponding Quot space.
\begin{thm}
With assumptions and notation as in the previous proposition, the moduli space $S_{E_0}^{\rm s}(P,\delta)=S_{E_0}(P,\delta)$ of $\delta$-stable pairs is isomorphic to the Quot space ${\rm Quot}^{\lambda,\#}(E_0,P)$ parametrizing quotients of $E_0$ with Hilbert polynomial $P$ in the heart $\coh^{\lambda,\#}(X)$. In particular, the Quot space is projective.
\end{thm}

\subsection{Framed sheaves}
We have a notion dual to pairs: {\em framed sheaves}. They are homomorphisms of the form $E\to E_0$.
\begin{defn}\label{framed-sheaf}
A {\em framed sheaf} is a coherent sheaf $E$ with Hilbert polynomial $P_E=P$, together with a nonzero {\em framing} $\alpha\colon E\to E_0$.
It is $\tau$-{\em stable} if $\ker \alpha$ is zero or pure of dimension $d$, and for every nonzero subsheaf $F\subset E$ of dimension $d
$,
\begin{enumerate}[label=(\roman*)]
    \item $p_F< p_E-\tau/r_E$ if $F\subset \ker \alpha$,
    \item $p_F-\tau/r_F< p_E-\tau/r_E$ otherwise.
\end{enumerate}
\end{defn}
\begin{thm}[Huybrechts-Lehn]
There is a projective coarse moduli space $F_{E_0}(P,\tau)$ of S-equivalence classes of $\tau$-semistable framed sheaves with Hilbert polynomial $P$. It contains an open subscheme $F_{E_0}^{\rm s}(P,\tau)$ as the fine moduli space of $\tau$-stable framed sheaves.
\end{thm}
When $\deg \tau\geqslant d$, the moduli space is isomorphic to a Quot scheme. Therefore, we again consider a small stability parameter $\tau$:
$$\deg \tau < \deg P =d.$$
Now, we let
\begin{equation*}
    \lambda=\frac{P-\tau}{r}.
\end{equation*}
In some cases, we can also identify stable framed sheaves as monomorphisms in a tilted heart.

\begin{prop}\label{framed-sheaf-injection}
 Suppose $E$ has Hilbert polynomial $P$ and $E_0\in \mathcal{T}^\lambda$.
 We also assume that $\tau$ is not a critical value. Then, $\alpha\colon E\to E_0$ is a monomorphism in $\cohltilt[1]$ if and only if $(E,\alpha\colon E\to E_0 )$ is a $\tau$-stable framed sheaf.
\end{prop}
The proof is similar to that of Proposition~\ref{pair-quotient}. For completeness, we include it here.

\begin{proof}
Suppose that $\alpha\colon E\to E_0$ is a monomorphism in $\cohltilt[1]$. By the assumption $E_0 \in \cT^{\lambda}$ and Lemma \ref{surjtil}, we know that $\alpha$ is a morphism in $\coh (X)$ with $E \in \cT^{\lambda}$ and $\ker (\alpha) \in \cF^{\lambda}$. Given a subsheaf $F\subset E$, $p_{E/F}>\lambda$, because $E\in \mathcal{T}^\lambda$. Then,
$p_F-\tau/r_F<\lambda$. If $F\subset \ker\alpha$, then $p_F\leqslant\lambda$, which is actually a strict inequality, since we assume $\tau$ is not critical. Therefore, $(E,\alpha\colon E\to E_0)$ is $\tau$-stable.

Conversely, we assume that $(E,\alpha\colon E\to E_0)$ is $\tau$-stable. First, $\tau$-stability implies that for any dimension $d$ quotient sheaf $Q$ of $E$, $p_Q-\tau/r_Q>\lambda$ or $p_Q>\lambda$. Therefore, $E\in \mathcal{T}^\lambda$. On the other hand, the $\tau$-stability also implies that if nonzero, $\ker (\alpha)$ is pure of dimension $d$ and has reduced Hilbert polynomial $\leqslant \lambda$. Furthermore, $\ker (\alpha)\in \mathcal{F}^\lambda$. Again by Lemma \ref{surjtil}, $\alpha$ is a monomorphism in $\cohltilt [1]$.
\end{proof}

\begin{rem}
This is a variant of \cite[Lemma 5.5]{rota2020quot}
.\end{rem}
Thus, by identifying two moduli functors, we can identify the moduli space of $\tau$-stable framed sheaves with the corresponding Quot space.
\begin{thm}
With assumptions and notation as in the previous proposition, the moduli space $F_{E_0}^{\rm s}(P,\tau)=F_{E_0}(P,\tau)$ of $\tau$-stable framed sheaves is isomorphic to the Quot space ${\rm Quot}^{\lambda, \#[1]}(E_0, P_{E_0}-P)$ parametrizing quotients of $E_0$ with Hilbert polynomial $P_{E_0}-P$ in the heart $\coh^{\lambda,\#}(X)[1]$. In particular, the Quot space is projective.
\end{thm}

\section{Change of Quot space under tilting}\label{wall-crossing}

In this section, we will prove a formula relating the moduli space of quotient husks and the Grothendieck's Quot scheme, which parameterize quotient sheaves supported in dimension no more than one. We follow Bridgeland's treatment of Hall algebra identities in \cite[\S 6]{Bridgeland11}.

\subsection{The stack of pairs}

We first introduce a modification of stack of sheaves with sections, which were constructed in \cite[\S 2.3]{Bridgeland11}. Let $X$ be a nonsingular projective variety over $\bC$ and $E_0 \in \Coh(X)$ fixed. We denote by $\cM$ the stack of coherent sheaves on $X$. It is an Artin stack, locally of finite type over $\bC$. There is another stack $\cM (E_0)$ with a morphism $q \colon \cM (E_0) \to \cM$ parameterizing pairs $(E, \alpha\colon E_0\to E)$. Indeed, the objects of $\cM (E_0)$ lying over a scheme $S$ are pairs $(E, \alpha)$ consisting of an $S$-flat coherent sheaf $E$ on $S \times X$ together with $\alpha \colon E_{0S} \to E$
where $E_{0S}$ denotes the pullback of $E_0$ under the projection $S \times X \to X$.
Let $f \colon T \to S$ be a morphism of schemes and $f_X = f \times \id_X$. Given an object $(F, \beta)$ lying over $T$, a morphism $\theta \colon (F, \beta) \to (E, \alpha)$ lying over $f$ is an isomorphism $\theta \colon f^{\ast}_X E \to F$ on $T \times X$ with $\theta \circ f^{\ast}_X \alpha = \beta \circ \kappa$, where the morphism $\kappa \colon f^{\ast}_X E_{0S} \to E_{0T}$ denotes the canonical isomorphism of pullbacks. The morphism $q$ of stacks  is defined by forgetting the data of the morphism $\alpha$ in the obvious way.

By an easy modification of the argument of \cite[Lemma 2.4]{Bridgeland11}, we have the following lemma.

\begin{lem}
The stack $\cM (E_0)$ is an Artin stack and the morphism $q$ is representable and of finite type.
\end{lem}

The following lemma is a result of the fibers of the morphism $q$.

\begin{lem} \label{ZarfriLem}
Let $E_0 \in \coh(X)$ be fixed. There is a stratification of $\cM$ by locally closed substacks $\cM_{r} \subset \cM$ such that the objects of $\cM (\bC)$ are $E \in \coh (X)$ with $\hom(E_0,E)=r$.
The pullback of $q \colon \cM (E_0) \to \cM$ to $\cM_r$ is a locally trivial fibration in the Zariski topology, with fiber $\bC^r$.
\end{lem}

\begin{proof}
Let $S$ be a scheme. Given an $S$-flat coherent sheaf $E$ on $S \times X$, we write $\mathfrak{h}om (E_{0S}, E)$ for the set-valued covariant funtor on $(S\mbox{-Schemes})^{\rm op}$, which associates to any $S$-scheme $f \colon T \to S$ the set $\Hom (f^{\ast}_X E_{0S}, f^{\ast}_X E)$ of $\cO_{T \times X}$-linear morphism . By a standard limit argument (cf.~\cite[(8.5.2), (8.8.2), (8.9.1), (11.2.6)]{EGAIV3}), we may assume that $S$ is noetherian. According to the results of Grothendieck (see \cite[Theorem 5.8]{Nitsure05} and references therein), there is a coherent sheaf $G(E_{0S}, E)$ on $S$ such that the funtor $\mathfrak{h}om (E_{0S}, E)$ is represented by the linear scheme
\begin{equation*} \label{linSch}
    \mathbf{Spec}\,(\sym_{\cO_S} G(E_{0S}, E)).
\end{equation*}
Then the remaining proof is essentially the same as in \cite[Lemma 2.5]{Bridgeland11}.
\end{proof}


\subsection{Motivic Hall algebra} \label{hull-subsec}


We are going to recall the notion of motivic Hall algebras. For a more detailed discussion,
we refer to \cite{Bridgeland11,Bridgeland12}.


We denote $\cC$ by the subcategory $\coh_{\leqslant 1} (X) $ of $\coh (X)$, and this corresponds to an open and closed substack $\cC \subset \cM$ by the usual abuse of notation. There exists a stack $\cC^{(2)}$ of short exact sequences in the category $\cC$. It comes with three distinguished morphisms $a_1, a_2$ and $b \colon \cC^{(2)} \to \cC$. These morphisms correspond to sending a short exact sequence $0 \to A_1 \to B \to A_2 \to 0$ to the sheaves $A_1, A_2$ and $B$ respectively. We remark that $(a_1, a_2)$ is of finite type \cite[Lemma 4.2]{Bridgeland12}.

The \emph{motivic Hall algebra}, denoted by $\rH (\cC)$, is the relative Grothendieck group $K (\St / \cC)$ over the stack $\cC$. By definition, it is defined to be the complex vector space spanned by isomorphism classes of symbols $[\cX \to \cC]$ where $\cX$ is an Artin stack of finite type over $\bC$ with affine geometric stabilizers
, modulo three relations: the scissor relations for finite disjoint satcks, geometric bijection
relations and Zariski fibration
relations (see \cite[Definition 3.10]{Bridgeland12}).

It is equipped with a \emph{noncommutative} product $\ast$ given explicitly by the rule
\begin{equation*}
    [\cX_1 \xrightarrow{f_1} \cC] \ast [\cX_2 \xrightarrow{f_2} \cC] = [\cZ \xrightarrow{b \circ h} \cC],
\end{equation*}
where $h$ is defined by the Cartesian diagram
\begin{equation*} \label{hullconv}
\begin{tikzcd}
\cZ \ar[r, "h"] \ar[d] & \cC^{(2)} \ar[r, "b"] \ar[d, "{(a_1, a_2)}"] & \cC. \\
\cX_1 \times \cX_2 \ar[r, "f_1 \times f_2"] & \cC \times \cC &
\end{tikzcd}
\end{equation*}
The unit is given by $1 = [\spec \bC \to \cC]$, which corresponds to the zero object in $\cC$ and
the product $\ast$ is associative \cite[Theorem4.3]{Bridgeland12}.

On the other hand, there is a natural grading on $\rH (\cC)$ by the monoid $\Delta$ consisting of classes of sheaves supported in dimension $\leqslant 1$. More precisely, let $N_1 (X)$ denote the abelian group of cycles of dimension one modulo numerical equivalence. We define the monoid by
\[
\Delta = \{ (\beta, n) \in N_1 (X) \oplus \bZ \mid \beta > 0 \text{ or } \beta = 0 \text{ and } n \geqslant 0 \}
\]
(cf.~\cite[\S 2.1]{Bridgeland11}). There are open and closed substacks $\cC_{\gamma} \subset \cC$, the stacks of objects of class $\gamma \in \Delta$. Thus elements of $\rH (\cC)$ are naturally graded by the monoid $\Delta$. An element $[f \colon \cX \to \cC]$ is {\em homogeneous} of degree $\gamma \in \Delta$ if $f$ factors through the substack $\cC_{\gamma}$.

\subsection{Laurent subsets}

Let us summarize sections 5.2 and 6.1 of \cite{Bridgeland11}. A subset $S \subset \Delta$ is \emph{Laurent} if for all $\beta \in N_1 (X)$, the collection $\{ n \in \bZ \mid (\beta, n) \in S\}$ is bounded below. Let $\Phi$ denote the set of all Laurent subsets.

For the $\Delta$-graded Hall algebra $\rH (\cC)$, we can use $\Phi$ to define a new algebra, denoted by $\rH (\cC)_{\Phi}$. Elements of such new algebra are of the form $a = \sum_{\gamma \in S} a_{\gamma}$
where $S \in \Phi$ and $a_{\gamma} \in \rH (\cC)_{\gamma} \subset \rH (\cC)$. 
There is a natural topology and product $\ast$ on $\rH (\cC)_{\Phi}$ induced by projection operators (see \cite[\S 5.2]{Bridgeland12}).

To define a stability condition, we fix an ample divisor $H$ on $X$. Given a class $\gamma = (\beta, n) \in \Delta$ define the slope $\mu (\gamma) = n (\beta \cdot H)^{- 1} \in (- \infty, \infty]$. If $\beta = 0$, we consider $\gamma$ to have slope $\infty$, and otherwise $\mu (\gamma) \in \bQ$.

Given an interval $I \subset (-\infty, \infty]$, define $\rSS (I) \subset \cC$  to be the full subcategory consisting of zero objects together with those one-dimensional sheaves whose Harder-Narasimhan factors all have slope in $I$ (see \cite[\S 6.1]{Bridgeland11}). We write $\rSS (I) = \rSS (\geqslant \mu)$ if $I = [\mu,  \infty]$. Then the following lemma follows from Lemmas 5.3, 6.2 and (31) of \cite{Bridgeland11}.

\begin{lemma} \label{inve-lem}
The subcategory $\rSS([\mu, \infty))$ defines an invertible element $1_{\rSS([\mu, \infty))}$ in $\rH (\cC)_{\Phi}$.
\end{lemma}




\subsection{Identities in the Laurent Hall algebra}

Let $\cT = \coh_{0} (X) = \rSS (\infty)$. Consider the torsion pair $(\cT, \cF \cap \cC)$ of $\cC$, where $\cF = \cT^{\perp}$. Then $\cC^{\#} = \left<\cF \cap \cC, \cT [-1]\right>$ is the tilt of $\cC$.

For the Grothendieck's Quot scheme $\Quot(E_0)$ and the moduli space of quotient husks (or limit stable pairs, Proposition \ref{husk-quotientCor}) $\QHusk(E_0)$, we introduce\footnote{The condition $E_0 \in \cF$ implies that if $\alpha \colon E_0 \to E$ is an epimorphism in $\coh (X)^{\#}$, then $E \in \cF \subset \coh (X)$ by Lemma \ref{surjtil}}
\begin{align*}
    \Quot(E_0)_{\leqslant 1} &\coloneqq \Quot(E_0) \cap \cC, \\
    \QHusk(E_0)_{\leqslant 1} &\coloneqq \QHusk(E_0) \cap \cC \text{ if } E_0 \in \cF
\end{align*}
which parameterize quotients of $E_0$ supported in dimension $\leqslant 1$. By the same argument of \cite[Lemma 2.6]{Bridgeland11}, we can view $\Quot(E_0)_{\leqslant 1}$ and $\QHusk(E_0)_{\leqslant 1}$ as open substacks of the moduli stack $\cC (E_0)$. Namely, these $\bC$-valued points are morphisms $E_0 \to E$ that are epimorphisms in the categories $\cC$ and $\cC^{\#}$ respectively. The morphisms
\begin{center}
    $\Quot(E_0)_{\leqslant 1} \to \cC$ and  $\QHusk(E_0)_{\leqslant 1} \to \cC$,
\end{center}
which are the restrictions to $\Quot(E_0)_{\leqslant 1}$ and $\QHusk(E_0)_{\leqslant 1}$ of $q \colon \cC (E_0) \to \cC$, define elements $\cQ_{\leqslant 1}$ and $\cQ_{\leqslant 1}^{\#}$ of $\cH (\cC)_{\Phi}$ by the similar argument of \cite[Lemma 5.5]{Bridgeland11}.

Given a substack $i \colon \cN \to \cC$, we write $1_{\cN} \coloneqq [\cN \to \cC]$ in $\rH (\cC)$.
Pulling back the morphism $q \colon \cC (E_0) \to \cC$ to $\cN \subset \cC$ given a stack $\cN (E_0)$ with a morphism $\cN (E_0) \to \cN$ and hence an element $1_{\cN}^{E_0} \coloneqq [\cN (E_0) \to \cC]$ in $\rH (\cC)$.
By abuse of notation, we use the same symbol for an open substack of $\cC$ and the corresponding full subcategory of $\cC$ defined by its $\bC$-valued points.

Following \cite{Bridgeland11} we establish the torsion pair and Quot 
space identities in the next two lemmas.

\begin{lem} \label{ideGB-L}
The following identities hold in the Laurent Hall algebra $\rH (\cC)_{\Phi}$.
\begin{enumerate}[label=(\alph*)]
    \item\label{ideGB-L1} $1_{\rSS(\geqslant \mu)} = 1_{\cT} \ast 1_{\rSS([\mu, \infty))}$.
    \item\label{ideGB-L2} $\lim_{\mu \to - \infty} \left(\cQ_{\leqslant 1} \ast 1_{\rSS(\geqslant \mu)} - 1_{\rSS(\geqslant \mu)}^{E_0} \right) = 0$.
    \item\label{ideGB-L3} $\lim_{\mu \to - \infty} \left(\cQ_{\leqslant 1}^{\#} \ast 1_{\rSS([\mu, \infty))} - 1_{\rSS([\mu, \infty))}^{E_0} \right) = 0$ if $E_0 \in \cF$.
\end{enumerate}
\end{lem}

The proof of Lemma \ref{ideGB-L} is essentially the same as in \cite[Proposition 6.5]{Bridgeland11}, noticing the boundedness of the Quot scheme and the moduli space of stable pairs. Remark that the geometric bijection relations plays an essential role, and we need the assumption $E_0 \in \cF$ of \ref{ideGB-L3} to use Lemma \ref{surjtil} instead of \cite[Lemma 2.3]{Bridgeland11}.

\begin{lem} \label{ideZF-L}
Assume that $E_0$ is locally free. There is an identity
\begin{center}
    $ 1_{\rSS(\geqslant \mu)}^{E_0} = 1_{\cT}^{E_0} \ast 1_{\rSS([\mu, \infty))}^{E_0}$ \quad in $\rH (\cC)_{\Phi}$.
\end{center}
\end{lem}

\begin{proof}
Form Cartesian squares
\begin{equation*}
\begin{tikzcd}
\cY \ar[d] \ar[r, "p"] & \cX \ar[r, "j"] \ar[d] & \cC^{(2)} \ar[r, "b"] \ar[d, "{(a_1, a_2)}"] & \cC \\
\cT (E_0) \times \rSS([\mu, \infty)) (E_0) \ar[r, "{(q, \id)}"] & \cT \times \rSS([\mu, \infty)) (E_0) \ar[r, "{(i, q)}"] & \cC \times \cC. &
\end{tikzcd}
\end{equation*}
Then $1_{\cT}^{E_0} \ast 1_{\rSS([\mu, \infty))}^{E_0}$ is represented by the composite morphism $b \circ j \circ p \colon \cY \to \cC$. Note that, by Lemma \ref{ZarfriLem}, the morphism of stacks $q \colon \cT (E_0) \to \cT$ is a Zariski fibration, with fiber over a sheaf $T$ being the vector space $\Hom (E_0, T)$. By pullback the same is true for the map $p$.

Since the morphism $(a_1, a_2)$ satisfies the iso-fibration property of \cite[Lemma A.1]{Bridgeland12}, the groupoid of $S$-valued points of $\cX$ is as follows. The objects are a short exact sequences of $S$-flat sheaves on $S \times X$
\begin{equation}\label{ZF-diag-L}
0 \to T \to E \xrightarrow{\gamma} F \to 0
\end{equation}
together with $T \in \cT$, $F \in \rSS([\mu, \infty))$ and a map $\alpha \colon E_{0S} \to F$.
The morphisms are isomorphisms of short exact sequences commuting with the
map $\alpha$.

Recall that by Lemma \ref{ideGB-L} \ref{ideGB-L1}, $b \circ s \colon \cZ \to \rSS (\geqslant \mu)$ induces an equivalence on $\bC$-valued points. On the other hand, consider a Cartesian diagram
\begin{equation*}
\begin{tikzcd}
\cW \ar[r, "h"] \ar[d]  & \rSS(\geqslant \mu) (E_0) \subset \cC (E_0) \ar[d, "q"] \\
\cZ \ar[r, "b \circ s"] & \rSS(\geqslant \mu) \subset \cC .
\end{tikzcd}
\end{equation*}
Since $b \circ s$ is a geometric bijection, so too is $h$. Thus the element $1_{\rSS(\geqslant \mu)}^{E_0}$ can be represented by the morphism $q \circ h$.

The groupoid of $S$-valued points of $\cW$ can be represented by the short exact sequences \eqref{ZF-diag-L} with a map  $\delta \colon E_{0S} \to E$.
Setting $\alpha = \gamma \circ \delta$ defines a morphism of stacks $\cW \to \cX$.  It is easy to see that this is a Zariski fibration, with fiber over a $\bC$-valued point of $\cX$ represented by a sequences \eqref{ZF-diag-L} with a map $\alpha$ being a vector space for $\Hom (E_0, T)$. Indeed, we have a long exact sequence
\begin{equation*}
    0 \to \Hom (E_0 , T) \to \Hom (E_0 , E) \to \Hom (E_0 , F) \to \Ext^1 (E_0, T)
\end{equation*}
on $X$. Since the dimension of support of $T$ is zero, so too is that of $E_0^{\vee} \otimes T$. For a locally free sheaf $E_0$, we get
\[
\Ext^1 (E_0, T) \cong H^1 (X, E_0^{\vee} \otimes T) = 0
\]
by the dimensional reason. Since $\cW \to \cX$ is the same fiber $\Hom (E_0 , T)$ as the map $p \colon \cY \to \cX$, the result follows from the Zariski relation $[\cW \to \cX \to \cC] = [\cY \to \cX \to \cC]$.
\end{proof}

We are now in a position to give the formula relating $\Quot(E_0)_{\leqslant 1}$ and $\QHusk(E_0)_{\leqslant 1}$.

\begin{thm}\label{App_thm}
Assume that $E_0 \in \cF$ and it is locally free. There is an identity
\begin{center}
    $\cQ_{\leqslant 1} \ast 1_{\cT} = 1_{\cT}^{E_0} \ast \cQ_{\leqslant 1}^{\#}$ \quad in $\rH (\cC)_{\Phi}$.
\end{center}
\end{thm}

\begin{proof}

By Lemma \ref{ideGB-L} \ref{ideGB-L1} and Lemma \ref{ideZF-L}, the expression \ref{ideGB-L2} of Lemma \ref{ideGB-L} can be rewritten
\begin{equation*}
    \cQ_{\leqslant 1} \ast 1_{\cT} \ast 1_{\rSS([\mu, \infty))} - 1_{\cT}^{E_0} \ast 1_{\rSS([\mu, \infty))}^{E_0} \to 0 \text{ as } \mu \to - \infty.
\end{equation*}
Multiplying \ref{ideGB-L3} of Lemma \ref{ideGB-L} on the left by $1_{\cT}^{E_0}$ gives
\begin{equation*}
    1_{\cT}^{E_0} \ast \cQ_{\leqslant 1}^{\#} \ast 1_{\rSS([\mu, \infty))} - 1_{\cT}^{E_0} \ast 1_{\rSS([\mu, \infty))}^{E_0} \to 0 \text{ as } \mu \to - \infty.
\end{equation*}
Thus
\begin{equation*}
    1_{\cT}^{E_0} \ast \cQ_{\leqslant 1}^{\#} \ast 1_{\rSS([\mu, \infty))} - \cQ_{\leqslant 1} \ast 1_{\cT} \ast 1_{\rSS([\mu, \infty))} \to 0 \text{ as } \mu \to - \infty.
\end{equation*}
By Lemma \ref{inve-lem}, we can multiply the inverse of $1_{\rSS([\mu, \infty))}$ and deduce the result.
\end{proof}

\begin{rem}
In \cite{toda20}, Toda studied the higher rank DT/PT  
correspondence, via stable objects in the derived category of coherent sheaves. He applied the integration map to the moduli stacks. By Behrend's result\cite{beh09}, the integrations are related to higher rank DT and PT invariants. The invariants are defined using the virtual fundamental classes, whose existence are guaranteed by the symmetric obstruction theories.

The moduli space of quotient husks/limit stable pairs can also be viewed as a version of higher rank PT moduli space. Over a Calabi--Yau 3-fold, we can also apply the integration map. However, the question whether the result is a deformation invariant remains, due to the absence of a result on a virtual fundamental class at the moment.
\end{rem}

{\it Acknowledgements.} We would like to thank Dingxin Zhang for helpful discussions. YL was supported by Grant 2017M620726 from China Postdoctoral Science Foundation. SSW thanks Shing-Tung Yau Center of Southeast University for providing a stimulating environment, and was supported by the Fundamental Research Funds for the Central Universities 2242020R10048.

\bibliography{main}{}
\bibliographystyle{alpha}
\end{document}